\documentclass[a4paper,12pt]{article}
\usepackage{amssymb,amsmath,amsthm,amscd,latexsym}

\newcommand{\Om} {\Omega}
\newcommand{\pa} {\partial}
\newcommand{\al} {\alpha}

\newcommand{\De} {\Delta}

\newcommand{\na} {\nabla}
\newcommand{\be} {\begin{equation}}
\newcommand{\ee} {\end{equation}}
\newcommand{\la} {\lambda}

\def\XXint#1#2#3{{\setbox0=\hbox{$#1{#2#3}{\int}$}
     \vcenter{\hbox{$#2#3$}}\kern-.5\wd0}}

\begin{document}
\title{ A note on the smoothness of energy-minimizing incompressible
deformations}

\author{{\Large Nirmalendu Chaudhuri} \\\\
Mathematical Sciences Research Institute\\
17 Gauss Way\\
Berkeley, CA 94720-5070, USA\\
chaudhur@msri.org\\
\&\\
Centre for Mathematics and its Applications\\
Australian National University\\
Canberra, ACT 0200\\
Australia\\
chaudhur@maths.anu.edu.au\\
}

\date{\empty}
\maketitle
\begin{center}{Abstract}
\medskip

{\footnotesize\rm \begin{tabular}{p{14.5cm}} In this note we prove
that any $W^{1,2}$ mapping $u$ in the plane that minimizes an
appropriate quasiconvex energy functional subject to the Jacobian
constraint ${\rm det}\,\na u=1$ a.e., are necessarily Lipschitz.
Furthermore we show that the minimizers corresponding to uniformly
convex energy are affine and give an example of non-affine
minimizers subject to affine boundary data corresponding to a
convex energy. We also discuss the regularity issues in dimension
greater than or equal to $3$.
\end{tabular} }
\end{center}

\section{Introduction}
Let $\Om\subset\Bbb R^2$ be a bounded incompressible material
body, that is, every $W^{1,2}$ deformation of $\Om$ locally
preserves its volume, in particular, the Jacobian of every such
deformation is 1 almost everywhere. For an incompressible
neo-Hookean material \cite{Og}, \cite{BOP} such as vulcanized
rubber, in the equilibrium, one is interested minimizing the
potential energy
\begin{equation}
I[u]:=\int_{\Om} F(\na u),
 \label{z1}
\end{equation}
for incompressible $W^{1,2}$ deformations $u:\Om\to\Bbb R^2$ with
prescribed boundary conditions corresponding to a given bulk
energy $F:\Bbb R^{2\times 2}\to \Bbb R$. The simplest $F$ is the
Dirichlet energy ${\displaystyle F(X):=\frac{1}{2}{\rm
tr}\,(X^TX)}$. Let us denote the incompressible or so-called the
area-preserving mappings
\begin{equation}
{\cal A}:=\{u\in W^{1,2}(\Om,\Bbb R^2)\,:\,\,{\rm det}\na u(x)=1,
{\rm a.\,e.}\,\,\,{\rm in}\,\,\Om\}, \label{z3}
\end{equation}
$u=(u^1,u^2)$, ${\displaystyle \na u=(u^{i}_{x_j})_{1\leq
i,j\leq2}}$, the gradient and ${\rm det}\,\na
u:=u^1_{x_1}u^2_{x_2}-u^1_{x_2}u^2_{x_1}$, is the Jacobian of $u$.
A function $f:\Bbb R^{m\times n}\to\Bbb R$ is said to be {\it
quasiconvex} if
$$
\int_{(0,1)^n}f(X+\na\phi(x))\,dx\geq f(X)
$$
for each $m\times n$ matrices $X$ and each smooth compactly
supported $\phi:(0,1)^n\to\Bbb R^m$. This definition was
introduced by Morrey \cite{Mo}, as a necessary and sufficient
condition for weak lower semicontinuity of the energy functional
associated to $f$ with respect to uniform convergence of Lipschitz
functions. Under the basic assumption that $F:\Bbb
R^{2\times2}\to\Bbb R$ is smooth, quasiconvex with quadratic
growth, together with the weak continuity \cite{Mu} of the
Jacobian, the functional $I$ admits {\it local} (see for example
{\cite{EG}) minimizers in the class ${\cal A}$. It remains a
difficult problem (due to the hard Jacobian constraint) to
understand the regularity properties of the local minimizes of
$I$. Under the additional assumption (the so-called uniform
quasiconvexity, see, \cite{Ev})
\begin{equation}
\int_{(0,1)^2}F(X+\na\phi(x))-F(X)\,dx\geq
C\int_{(0,1)^2}|\na\phi|^2dx \label{z2}
\end{equation}
 for some $C>0$, for each $2\times 2$ matrices $X$ and each smooth
 compactly supported $\phi:(0,1)^2\to\Bbb R^2$, Evans and Gariepy
\cite{EG}
 proved that any {\it non-degenerate, Lipschitz area-preserving} local
minimizers of $I$ are $C^{1,\al}$ on a dense open subset. It
remains to understand whether area-preserving local minimizers are
Lipschitz. Here we consider only the {\it global minimizers}. A
map $u\in{\cal A}$ is said to be {\it global minimizer} of $I$
subject to its own boundary  if
$$
I[u]\leq I[v],\,\,\,\,{\rm for}\,\,{\rm all}\,\,\,v\in{\cal A}.
$$
\medskip

\noindent
 {\bf Theorem 1.1.} {\it Let $F:\Bbb R^{2\times 2}\to
\Bbb R$ be uniformly quasiconvex, $C^2$, and $D^2 F$ is bounded.
Then global minimizers of $I$ in the area-preserving class ${\cal
A}$ are Lipschitz. Furthermore, if $F$ is uniformly convex and
frame indifference then the minimizers are affine}
\medskip

The proof follows by reducing the minimization problem to a
partial differential relation of the form
$$
\na u(x)\in K\,\,\,\,{\rm a.e.}\,\,\,\,\,\,{\rm in}\,\,\Om,
$$
for suitable subset $K$ of $\Bbb R^{2\times 2}$. As a consequence
of this observation and a theorem of M\"uller and \v{S}ver\'{a}k
\cite{MS} on convex integration (also see, \cite{DM}), we give an
example of non-affine minimizers subject to affine boundary data
corresponding to a convex (non-uniform) energy.

\section{Proof of Theorem 1.1.}
We recall the set of area-preserving mappings
$$
{\cal A}:=\{u\in W^{1,2}(\Om,\Bbb R^2)\,:\,\,{\rm det}\na u(x)=1,
{\rm a.\,e.}\,\,\,{\rm in}\,\,\Om\}.
$$
Let $SL(2):=\{P\in\Bbb R^{2\times2}\,:\,{\rm det}\,P=1\}$, the
special linear group and
\begin{equation}
Z_{\min}(F):=\{Q\in SL(2)\,:\,F(Q)=\min_{P\in SL(2)}F(P)\},
\label{z4}
\end{equation}
be the minimizing set. Since $F$ is uniformly quasiconvex, the
minimizing set $Z_{\min}(F)$ is non-empty and compact. Without
loss of generality assume $|\Om|=1$. Observe that for any $v\in
{\cal A}$, we have
\begin{equation}
I[v]=\int_{\Om}F(\na v)\geq \min\{F(P)\,:\,{\rm det}\,P=1\}.
\label{z5}
\end{equation}
Let $Q\in Z_{\min}(F)$ and $u(x)=Qx$, be a linear deformation.
Then $u\in {\cal A}$ and
$$
I[u]=F(Q)=\min_{P\in SL(2)} F(P).
$$
Therefore
\begin{equation}
\min_{v\in{\cal A}} I[v]=\min_{P\in SL(2)}F(P) \label{z8}.
\end{equation}
Hence $u\in {\cal A}$ is a minimizer of $I$ if and only if it
satisfies the partial differential inclusion
\begin{equation}
\na u(x)\in Z_{\min}(F) \,\,\,\,{\rm a.e.}\,\,\,\,\,\,{\rm
in}\,\,\Om, \label{z6}
\end{equation}
Since $Z_{\min}(F)$ is compact, $u$ is Lipschitz. This proves
first part of the theorem.
\bigskip

\noindent
 {\bf Lemma 2.1.} {\it Let $F:\Bbb R^{2\times2}\to\Bbb R$
be uniformly convex, that is, $D^2F(X)Y:Y\geq 2\la |Y|^2$ for some
$\la>0$, for all $X,Y\in\Bbb R^{2\times2}$. Suppose further,
$F(RX)=F(X)$ for each $X \in\Bbb R^{2\times2}$, and each rotations
$R$. Then the minimizing set $Z_{\min}(F)$ is simply the coset
$SO(2)P$, for some ${\rm det}\,P=1$, where $SO(2):=\{R\in\Bbb
R^{2\times2}\,:\,R^TR=Id,\,{\rm det}\,R=1\}$, the special
orthogonal group.}
\medskip

{\it Proof.} Here we follow the standard uniqueness arguments.
Suppose there exists $Q_1, Q_2\in SL(2)$ that
$F(Q_1)=F(Q_2)=\min_{SL(2)}F(P)$. Since $F$ is frame indifferent,
it follows that any
$$Q\in K:=SO(2)Q_1\,\cup\,SO(2)Q_2$$
also minimizes $F$ over $SL(2)$. We claim that $Q_1$ and $Q_2$ are
conformally equivalent, i.e., $Q_1=RQ_2$, for some $R\in SO(2)$.
Suppose, $Q_1$ and $Q_2$ are not conformally equivalent. Since
${\rm det}\, Q_1={\rm det}\, Q_2=1$, a simple calculation shows
that the cosets $SO(2)Q_1$ and $SO(2)Q_2$ are rank-one connected.
Therefore for each $P_1\in SO(2)$ there exists $P_2\in SO(2)$ such
that
$$
{\rm det}\,(P_1Q_1-P_2Q_2)=0.
$$
Since $X\mapsto {\rm det}\,X$ is linear along rank-one directions,
it follows that
$$
{\rm det}\,(\la P_1Q_1+(1-\la)P_2Q_2)
=1,\,\,\,{\rm for}\,\,\,{\rm all}\,\,0\leq\la\leq1.
$$
From the uniform convexity of $F$, we have
\begin{equation}
F(X)\geq F(Y)+DF(Y):(X-Y) +\la|X-Y|^2 ,\,\,\,{\rm for}\,\,\,{\rm
all}\,\,X,\, Y\,\in\Bbb R^{2\times 2}, \label{z7}
\end{equation}
 for some $\la>0$, where $A : B :={\rm tr}(A^TB)$ is the scalar
product. By taking $X=P_1Q_1$ and $Y=(P_1Q_1+P_2Q_2)/2$, and vice
versa and adding these two inequalities, we obtain
$$
F(P_1Q_1)+F(P_2Q_2)\geq 2F\left(\frac{P_1Q_1+P_2Q_2}{2}\right)
+2\la|P_1Q_1-P_2Q_2|^2.
$$
Since $F(P_1Q_1)=F(P_2Q_2)=\min_{SL(2)} F$ and
$(P_1Q_1+P_2Q_2)/2\in SL(2)$, it follows that $P_1Q_1=P_2Q_2$, a
contradiction. Hence the minimizing set $Z_{\min}(F):=\{P\in
 SL(2)\,:\,F(P)=\min_{SL(2)}F(Q)\}$ is just one copy of the special
orthogonal group $SO(2)$. This proves the lemma.\qed
\medskip

Suppose $F$ is uniformly convex and $u\in{\cal A}$ is a minimizer
of $I$. Then by (\ref{z6}) and lemma 2.1, it follows that
\begin{equation}
\na u(x)\in SO(2)P \,\,\,\,{\rm a.e.}\,\,\,\,\,\,{\rm in}\,\,\Om,
\label{z9}
\end{equation}
for some ${\rm det}\,P=1$. From the Liouville Theorem of
Reshetnyak \cite{Re}, it follows that $u$ is affine. However for
the convenience of the readers we give a proof, which is due to
Kinderlehrer \cite{Ki} (the same proof works in all dimensions).
Let us make the change of variables, $v:P^{-1}(\Om)\to\Bbb R^2$,
by $v(P^{-1}x)=u(x)$, $x\in\Om$. Then $\na v(P^{-1}x)=\na
u(x)\,P$. Hence
$$
\na v(y)\in SO(2) \,\,\,\,{\rm a.e.}\,\,\,\,\,\,{\rm
in}\,\,\,P^{-1}(\Om).
$$
Since ${\rm cof}\,Q=Q$ on $SO(2)$, and
$$
{\rm div}\,{\rm cof}\,\na v=0,
$$
({\rm div} is taken in each rows) it follows that $v$ is harmonic,
i.e., $\De v=(\De v^1, \De v^2)=(0,0)$ and hence smooth. Since
$|\na v|^2=2$, the identity
$$
\frac{1}{2}\De|\na v|^2=|\na^2v|^2+\na v: \na\De v
$$
yields $\na^2 v=0$ in $P^{-1}(\Om)$ and hence $u$ is affine in
$\Om$. This proves the theorem.\qed
\bigskip

\noindent {\bf Remark 1.} The proof shows that the energy
minimizing volume-preserving $W^{1,n}$ deformations on bounded
open subsets of $R^n$, for $n\geq 2$ are Lipschitz.
\medskip

\noindent {\bf Remark 2.}  However, for $n\geq 3$, we are unable
to conclude whether minimizers corresponding to frame indifferent
uniformly convex functions on $\Bbb R^{n\times n}$ are necessarily
affine. Let us briefly discuss the case $n=3$. For any given
$3\times 3$ matrices $Q_1$ and $Q_2$ with determinant $1$, the
cosets $SO(3)\,Q_1$ and $SO(3)\,Q_2$ are not necessarily rank-one
connected (this is the main difference with the two dimension), so
the above proof fails to conclude $SO(3)\,Q_1=SO(3)\,Q_2$. If the
cosets are not rank-one connected (for example, the $SO(3)$ and
$SO(3)\,Q_2$, $Q_2={\rm diag}\,(\la_1,\la_2,\la_3)$,
$0<\la_1\leq\la_2\leq\la_3$, $\la_2\neq 1$, $\la_1\la_2\la_3=1$
are not rank-one connected) it is natural to determine whether any
Lipschitz map $u:\Om\subset\Bbb R^3\to\Bbb R^3$ satisfying
\begin{equation}
 \na u(x) \in Z:=SO(3)\,Q_1\,\cup\,SO(3)\,Q_2\,\,\,\,{\rm
a.e.}\,\,\,\,\,\,{\rm in}\,\,\,\Om\subset\Bbb R^3 \label{z11}
\end{equation}
are necessarily affine. From the separation lemma \cite[Lemma
2.4]{CM} it is enough to show that such solutions $u$ are
$W^{2,2}_{\rm loc}$. In order to obtain such regularity one
usually tries to find a suitable system of partial differential
equations for $u$ satisfying (\ref{z11}). However, it follows that
there are no uniformly elliptic system of PDEs (as the set $Z$ is
not {\it strongly incompatible}, see, \cite{CM}) for $u$ in
(\ref{z11}). This suggests that there are no obvious way of
getting $W^{2,2}$ regularity.

\section{Non-affine Minimizers}

In this section we show that there are convex functions $F$ for
which the functional $I$ in (\ref{z1}) admits area-preserving
non-affine minimizers even with prescribed affine boundary data.
To obtain such minimizers, idea is to look for a smooth convex
(not uniform) function $F$ such that the minimizing set
$Z_{\min}(F)$ strictly contains two copies of $SO(2)$, which are
rank-one connected. Then trivially one obtains non-affine
minimizers, for example, simple laminates. But interestingly, by a
theorem of M\"uller and \v{S}ver\'{a}k \cite{MS}, on convex
integration, we can find non-affine minimizers with prescribed
affine boundary.
\medskip

Let $H:={\rm diag}\,(\la,\mu)$, be a diagonal matrix such that
$0<\la<1<\mu$ and $\la\mu=1$. Set
$$K:=SO(2)\,\cup\,SO(2)\,H,
$$
the two wells. Then a simple calculation shows that each matrix in
$K$ is rank-one connected with exactly two other matrices in $K$.
Define $F:\Bbb R^{2\times2}\to\Bbb R$ by
$$
F:=\sup\{g\,:\,g\,\,{\rm convex}\,\,{\rm on}\,\,\Bbb R^{2\times
2},\,\,\,\, g\leq {\rm dist}^2(\cdot\,,\,K)\},
$$
the convex envelope of the square of the distance function ${\rm
dist}(\cdot\,,\,K)$. Therefore $F$ is smooth, convex and the
second derivative of $F$ uniformly bounded. Notice that $F(X)=0$
if and only if
$$
X\in K^{{\rm c}}=\left\{ \left(\begin{array}{ccc}
x_1 & -x_2 \\
x_2 & x_1 \\
\end{array} \right) +
\left(\begin{array}{ccc}
y_1 & -y_2 \\
y_2 & y_1 \\
\end{array} \right)
\left(\begin{array}{ccc}
\la & 0 \\
0& \mu\\
\end{array} \right)
: x, y\in\Bbb R^2,\,\,|x|+|y|\leq 1 \right\},
$$
the convex hull of the set $K$. Therefore the minimizing set
$Z_{\min}(F)=K^c\cap\,SL(2)$ given by
$$
Z_{\min}(F)=\left\{ \left(\begin{array}{ccc}
x_1+\la y_1 & -x_2-\mu y_2 \\
x_2+\la y_2& x_1+\mu y_1 \\
\end{array} \right)
: |x|+|y|\leq 1,\, |x|^2+|y|^2+(\la+\mu)\langle x, y\rangle=1
\right\},
$$
is the so-called rank-one convex hull (see \cite{Sv} for more
details about rank-one convex or quasiconvex hulls of general two
wells energy) of $K$. It is clear that the set $K$ is strictly
contained in $Z_{\min}(F)$. For $R\in Z_{\min}(F)\setminus K$.
M\"uller and \v{S}ver\'{a}k \cite{MS} on convex integration, the
following boundary value partial differential inclusion
\begin{equation}
\left\{\begin{array}{ccc}\na u(x)\in SO(2)\,\cup\,SO(2)\,H
\,\,\,\,\,\,{\rm a.e.}\,\,\,\,\,\,{\rm
in}\,\,\,\Om\\\\
u(x)=Rx+b\,\,\,\,\,\,\,\,{\rm on}\,\,\,{\pa\Om}
\end{array}\right.
\label{z10}
\end{equation}
admit solutions.
Therefore solutions to the problem (\ref{z10}) are clearly
non-affine and minimizes of the energy functional $I$ in
{\ref{z1}) over the class of functions
$${\cal A}_R:=\{u\in
W^{1,2}(\Om,\Bbb R^2)\,:\,{\rm det}\na u(x)=1, {\rm
a.\,e.}\,\,\,{\rm in}\,\,\Om,\,\,u(x)=Rx+b\,\,\,{\rm
on}\,\,\partial\Om\}.
$$
This shows that the uniform convexity assumption in the Theorem
1.1 is sharp.




\end{document}